\documentclass[]{aspm}

\articleinfo{}{}{}

\usepackage{verbatim}
\usepackage{amssymb}
\usepackage{amsbsy}
\usepackage{amscd}
\usepackage{amsmath}
\usepackage{amsthm}
\usepackage[mathscr]{eucal}

\usepackage{latexsym}
\usepackage{url}

\newtheorem{theorem}{Theorem}
\newtheorem{remark}{Remark}

\numberwithin{defn}{section}

\def\R{{\bf R}}

\def\N{{\bf N}}

\def\d{\displaystyle}
\def\e{{\varepsilon}}

\def\wt{\widetilde}

\def\p{\partial}

\newcommand{\LR}[1]{{\langle {#1} \rangle }}

\title[]{Recent developments
on the lifespan estimate for classical solutions
of nonlinear wave equations\\
in one space dimension}

\dedicatory{\footnotesize \it Dedicated to Professor Tohru  Ozawa
on the occasion of his sixties birthday}

\author[]{Hiroyuki Takamura}

\address{Mathematical Institute,
Tohoku University,
Aoba, Sendai 980-8578, Japan.}

\email{hiroyuki.takamura.a1@tohoku.ac.jp}

\rcvdate{}
\rvsdate{}



\begin{document}

\begin{abstract}
In this paper, we overview the recent progresses on the lifespan estimates
of classical solutions of the initial value problems
for nonlinear wave equations in one space dimension.
There are mainly two directions of the developments on the model equations
which ensure the optimality of the general theory.
One is on the so-called \lq\lq combined effect" of two kinds of the different nonlinear terms,
which shows the possibility to improve the general theory.
Another is on the extension to the non-autonomous nonlinear terms
which includes the application to nonlinear damped wave equations
with the time-dependent critical case.
\end{abstract}

\maketitle

\section{Introduction.}
In order to illustrate our purpose, let us turn back to the general theory
for nonlinear wave equations in one space dimension
which was introduced by Li, Yu and Zhou \cite{LYZ91,LYZ92} more than 30 years ago.
\par
We consider the initial value problem of the form;
\begin{equation}
\label{GIVP}
\left\{
\begin{array}{ll}
\d  u_{tt}- u_{xx}=H(u,u_t,u_x,u_{xx},u_{xt}) &\quad \mbox{in}\ \R\times (0,T),\\
 u(x,0)= \e f(x) , \quad u_t(x,0)= \e g(x) &\quad \mbox{for}\ x\in\R,
\end{array}
\right.
\end{equation}
where $T>0$, $f,g\in C_0^{\infty}(\R)$ and $\e>0$ is a sufficiantly small parameter. 
Let $\wt{\lambda}=(\lambda;\ (\lambda_i), i=0,1; \ (\lambda_{ij}), i,j,=0,1, i+j\ge1).$
Assume that $H=H(\wt{\lambda})$ is a sufficiently smooth function with 
\[
H(\wt{\lambda})=O(|\wt{\lambda}|^{1+\alpha})
\]
in a neighborhood of $\wt{\lambda}=0$, where $\alpha \in \N$.
Let us define the lifespan $\wt{T}(\e)$ as the maximal existence time of the classical solution of (\ref{GIVP}) with arbitrary fixed data.
We are interested in the long-time stability of the trivial solution due to the fact that
we cannot expect any time-decay of the solution of the free wave equation in one space dimension.
Indeed, the general theory is to express the lower bound of $\wt{T}(\e)$ by means of the smallness
of the initial data, i.e. $\e$, for which, Li, Yu and Zhou \cite{LYZ91, LYZ92} obtained
\begin{equation}
\label{lifespan_general}
\wt{T}(\e)\ge 
\left\{
\begin{array}{lll}
C\e^{-\alpha/2} & \mbox{in general case,}\\
C\e^{-\alpha(\alpha+1)/(\alpha+2)} & \mbox{if}\ \d \int_{\R}g(x)dx = 0,\\
C\e^{-\min\{\beta_0/2,\alpha\}} & \mbox{if}\ \p_{u}^{\beta}H(0)=0\ \mbox{for}\ 1+\alpha\le \forall \beta \le\beta_0,
\end{array}
\right.
\end{equation}
where $C$ is a positive constant independent of $\e$.
This result has been expected complete more than 30 years.
\par
Beyond the general theory, our interest went to its optimality or to extending the general theory
by studying the morel problem;
\begin{equation}
\label{IVP_gcombined}
\left\{
\begin{array}{ll}
	\d u_{tt}-u_{xx}=A(x,t)|u_t|^p|u|^q+B(x,t)|u|^r
	&\mbox{in}\quad \R\times(0,T),\\
	u(x,0)=\e f(x),\ u_t(x,0)=\e g(x),
	& x\in\R,
\end{array}
\right.
\end{equation}
where $p,q,r>1$ ($q$ could be zero) and $A,B$ are non-negative functions of space-time variables.
Let us write the lifespan of classical solutions of (\ref{IVP_gcombined}) by $T(\e)$.
According to the series of our studies on the estimates of $T(\e)$,
we will see later that the general theory can be improved in some case,
which arises from the constant coefficient case.
This part is presented in Section 2.
Moreover, the principle in extending the nonlinear term $H$ in (\ref{GIVP}) of the general theory
to the non-autonomous one $H=H(x,t,u,u_t,u_x,u_{xx},u_{xt}) $
must be initiated by variable coefficient case.
This part is presented in Section 3.


\section{Case of constant coefficients and the combined effect.} 
In this section, we assume that
\[
A(x,t)\equiv A_0\quad\mbox{and}\quad B(x,t)\equiv B_0.
\]
where $A_0$ and $B_0$ are non-negative constants.

\subsection{The generalized combined effect.}
When $A_0=0$ and $B_0>0$,
Zhou \cite{Zhou92} obtained the estimates of $T(\e)$;
\begin{equation}
\label{lifespan_A=0}
T(\e)\sim
\left\{
\begin{array}{lll}
C\e^{-(r-1)/2} &\mbox{if}&\d \int_{\R} g(x)dx \neq0,\\
C\e^{-r(r-1)/(r+1)} &\mbox{if}&\d \int_{\R} g(x)dx =0.
\end{array}
\right.
\end{equation}
Here we denote the fact that there are positive constants,
$C_1$ and $C_2$, independent of $\e$ satisfying $E(\e,C_1)\le T(\e)\le E(\e,C_2)$
by $T(\e)\sim E(\e,C)$.
The classification by total integral of the initial speed $g$ is caused by strong Huygens' principle
such as (\ref{Huygens}).
On the other hand, when $A_0>0$ and $B_0=0$, we have
\begin{equation}
\label{lifespan_B=0}
T(\e)\sim C\e^{-(p+q-1)}
\end{equation}
For $q=0$, the upper bound in this estimate was obtained by Zhou \cite{Zhou01},
and the lower bound is due to Kitamura, Morisawa and Takamura \cite{KMT23}.
For $q>1$, (\ref{lifespan_B=0}) was verified by Zhou \cite{Zhou} for the upper bound
with integer $p,q$ satisfying $p\ge1,q\ge0,p+q\ge2$,
and by Li,Yu and Zhou \cite{LYZ91,LYZ92}
for the lower bound with integer $p,q$ satisfying $p+q\ge2$
including more general but smooth terms.
Note that \cite{Zhou} is a preprint version of Zhou \cite{Zhou01}
in which only the case of $q=0$ is considered.
But it is easy to apply its argument to the case of $q>1$.
The lower bound in this case is due to Kido, Sasaki, Takamatsu and Takamura \cite{KSTT}.
\par
Therefore the natural expectation is that
\begin{equation}
\label{lifespan_conjecture}
T(\e)\sim
\left\{
\begin{array}{ll}
C\e^{-\min\{(p+q-1),(r-1)/2\}} & \mbox{if}\ \d\int_{\R}g(x)dx\not=0,\\
C\e^{-\min\{(p+q-1),r(r-1)/(r+1)\}} & \mbox{if}\ \d\int_{\R}g(x)dx=0
\end{array}
\right.
\end{equation}
in the case where $A_0>0$ and $B_0>0$.
But, surprisingly, we have the following fact.

\begin{theorem}[Morisawa, Sasaki and Takamura \cite{MST,MST_erratum}, Kido, Sasaki, Takamatsu and Takamura \cite{KSTT}] 
\label{thm:gcombined}
The conjecture (\ref{lifespan_conjecture}) is true except for the case where
\begin{equation}
\label{gcombined}
\int_{\R}g(x)dx=0\quad\mbox{and}\quad \frac{r+1}{2}<p+q<r.
\end{equation}
In this case, we have that
\begin{equation}
\label{lifespan_gcombined}
T(\e)\sim C\e^{-(p+q)(r-1)/(r+1)},
\end{equation}
which is strictly shorter than the second case in (\ref{lifespan_conjecture}).
\end{theorem}

\par
The brief proof of (\ref{lifespan_gcombined}) will appear at the end of this section.
We shall call this special phenomenon by \lq\lq generalized combined effect" of two nonlinearities.
The original \lq\lq combined effect", which means the case of $q=0$, was first observed by
Han and Zhou \cite{HZ14} which targeted to show the optimality of the result
by Katayama \cite{Katayama01} on the lower bound of the lifespan of classical solutions
of nonlinear wave equations with a nonlinear term $u_t^3+u^4$ in two space dimensions
including more general nonlinear terms.
It is known that $T(\e)\sim\exp\left(C\e^{-2}\right)$ for the nonlinear term
$u_t^3$ and $T(\e)=\infty$ for the nonlinear term $u^4$,
but Katayama \cite{Katayama01} obtained only a much worse estimate
than their minimum  as $T(\e)\ge c\e^{-18}$.
Surprisingly, more than ten years later, Han and Zhou \cite{HZ14} showed that
this result is optimal as $T(\e)\le C\e^{-18}$.
They also considered (\ref{IVP_gcombined}) with $q=0$ for all space dimensions $n$ bigger than 1
and obtain the upper bound of the lifespan.
Its counter part, the lower bound of the lifespan, was obtained by
Hidano, Wang and Yokoyama \cite{HWY16} for $n=2,3$.
See the introduction of \cite{HWY16} for the precise results and references.  
We note that the estimate (\ref{lifespan_gcombined}) with $q=0$ coincides with the lifespan estimate 
for the combined effect in \cite{HZ14, HWY16} if one sets $n=1$ formally.
Indeed, \cite{HZ14} and \cite{HWY16} showed that
\begin{equation}
\label{lifespan_highD}
T(\e)\sim C\e^{-2p(r-1)/\{2(r+1)-(n-1)p(r-1)\}}
\end{equation}
holds for $n=2,3$ provided
 \begin{equation}
 \label{condition_combined}
(r-1)\{(n-1)p-2\}<4,\ 2\le p\le r\le 2p-1,\ r>\frac{2}{n-1}.
\end{equation}
Later, Dai, Fang and Wang \cite{DFW19} improved the lower bound of lifespan
for the critical case in \cite{HWY16}.
They also show that $T(\e)<\infty$ for all $p,r>1$ in case of $n=1$, i.e. (\ref{IVP_gcombined})
with $q=0$.
For the non-Euclidean setting of the results above, see Liu and Wang \cite{LW20} for example,
in which the application to semilinear damped wave equations is included.

\subsection{Comparison with the general theory.}
Here we strongly remark that our estimate in (\ref{lifespan_gcombined}) is better than that of the general theory by Li, Yu and Zhou \cite{LYZ91, LYZ92}
in the case of (\ref{gcombined}) with integer $p,q,r\ge2$.
Because our result on the lower bound of the lifespan can be established
also for the smooth terms as $u_{tt}-u_{xx}=u_t^pu^q+u^r$.
We note that there are infinitely many examples of $(p,q,r)=(m,m,2m+1)$ as the inequality
\[
\frac{r+1}{2}=m+1<p+q=2m<r=2m+1
\]
holds for $m=2,3,4,\ldots$.
This fact shows a possibility to improve the general theory.
We also note that, even for the original combined effect of $q=0$,
the integer points satisfying (\ref{condition_combined}) are
$(p,r)=(2,3),(3,3),(3,4)$ for $n=2$ and $(p,r)=(2,2)$ for $n=3$,
but (\ref{lifespan_highD}) with $p=r$ agrees with the case of $A_0=0$ and $B_0>0$.
See Introduction of Imai, Kato, Takamura and Wakasa \cite{IKTW}
for references on the case of $A_0=0$ and $B_0>0$.
Hence one can say that only the lifespan estimates with $(p,r)=(2,3),(3,4)$ for $n=2$
are essentially in the combined effect case.
If $q\neq0$, $p$ is replaced with $p+q$ in the results above.
Therefore it has less meaningful to consider (\ref{IVP_gcombined})
in higher space dimensions, $n\ge2$, if we discuss the optimality of the general theory.
In spite of this situation, Han and Zhou \cite{HZ14} studied
\[
u_{tt}-\Delta u=u|u_t|^{p-1}+u|u|^{q-1}
\quad\mbox{in}\ \R^2\times(0,T)
\]
with $1<p\le4,q>1$ to show the blow-up part of the generalized combined effect.
\par
Of course, some special structure of the nonlinear terms
such as \lq\lq null condition" guarantees the global-in-time existence.
See Tartar \cite{Tartar}, Bianchini and G. Staffilani \cite{BS},
Nakamura \cite{Nakamura14}, Luli, Yang and Yu \cite{LYY18},
Zha \cite{Zha20, Zha22} for this direction.
But we are interested in the optimality of the general theory.

\par
From now on, let us make sure the fact above.
If one applies the result of the general theory (\ref{lifespan_general}) to our problem (\ref{IVP_gcombined}) with
\begin{equation}
\label{H_special}
H(u,u_t,u_x,u_{xx},u_{xt})=u_t^pu^q+u^r\quad\mbox{with}\  p,q,r\in\N,
\end{equation}
one has the following estimates in each cases.
\begin{itemize}
\item
When $p+q<r$,
\par
then, we have to set $\alpha=p+q-1$ and $\beta_0=r-1$ which yield that
\[
\wt{T}(\e)\ge
\left\{
\begin{array}{ll}
C\e^{-(p+q-1)/2} & \mbox{in general},\\
C\e^{-(p+q)(p+q-1)/(p+q+1)} & \mbox{if}\ \d\int_{\R}g(x)dx=0,\\
C\e^{-\min\{(r-1)/2,p+q-1\}} &
\begin{array}{l}
\mbox{if $\p_u^\beta H(0)=0$}\\
\mbox{for $p+q\le\forall\beta\le r-1$}.
\end{array}
\end{array}
\right.
\]
We note that the third case is available for (\ref{H_special}).
Therefore, for $p+q\le(r+1)/2$, we obtain that
\[
\wt{T}(\e)\ge c\e^{-(p+q-1)}
\] 
whatever the value of $\d\int_{\R}g(x)dx$ is.
On the other hand, for $(r+1)/2<p+q$, i.e. 
\[
\frac{r-1}{2}<p+q-1,
\]
we obtain
\[
\wt{T}(\e)\ge
\left\{
\begin{array}{ll}
C\e^{-(r-1)/2} & \mbox{if}\ \d\int_{\R}g(x)dx\neq0,\\
C\e^{-\max\{(r-1)/2,(p+q)(p+q-1)/(p+q+1)\}} & \mbox{if}\ \d\int_{\R}g(x)dx=0.\\
\end{array}
\right.
\]
\item
When $p+q\ge r$,
\par
then, similarly to the case above, we have to set $\alpha=r-1$, which yields that
\[
\wt{T}(\e)\ge
\left\{
\begin{array}{ll}
C\e^{-(r-1)/2} & \mbox{in general},\\
C\e^{-r(r-1)/(r+1)} & \mbox{if}\ \d\int_{\R}g(x)dx=0,\\
C\e^{-\min\{\beta_0/2,(r-1)\}} & \mbox{if $\p_u^\beta H(0)=0$ for $r\le\forall\beta\le\beta_0$}.
\end{array}
\right.
\]
We note that the third case does not hold for (\ref{H_special}) by $\p_u^r H(0)\neq0$.
\end{itemize}
As a conclusion, for the special nonlinear term in (\ref{H_special}),
the result of the general theory is
\[
\wt{T}(\e)\ge
\left\{
\begin{array}{ll}
C\e^{-(p+q-1)} & \mbox{for}\ p+q\le\d\frac{r+1}{2},\\
C\e^{-(r-1)/2} & \mbox{for}\ \d\frac{r+1}{2}\le p+q
\end{array}
\right.
\mbox{if}\ \int_{\R}g(x)dx\not=0\\
\]
and
\[
\wt{T}(\e)\ge
\left\{
\begin{array}{l}
C\e^{-(p+q-1)}\\
\qquad\d\mbox{for}\ p+q\le\d\frac{r+1}{2},\\
C\e^{-\max\{(r-1)/2,(p+q)(p+q-1)/(p+q+1)\}}\\
\qquad\d\mbox{for}\ \d\frac{r+1}{2}\le p+q\le r,\\
C\e^{-r(r-1)/(r+1)}\\
\qquad\mbox{for}\ r\le p+q\\
\end{array}
\right.
\mbox{if}\ \d\int_{\R}g(x)dx=0.
\]
Therefore a part of our results in (\ref{lifespan_gcombined})
is larger than the lower bound of $\wt{T}(\e)$ in that case.
If one follows the proof of Theorem \ref{thm:gcombined}, one can find that it is easy to see that
our results on the lower bounds also hold for a special term (\ref{H_special})
by estimating the difference of nonlinear terms from above
after employing the mean value theorem.
This fact indicates that we still have a possibility to improve the general theory
in the sense that the optimal results in (\ref{lifespan_gcombined}) should be included at least.

\subsection{An improvement of the general theory.}
Inspired by Theorem \ref{thm:gcombined}, Takamatsu \cite{Takamatsu} has
recently proved the following theorem which is an improvement of the general theory
related to the generalized combined effect.

\begin{theorem} [Takamatsu \cite{Takamatsu}]
\label{thm:Takamatsu}
The fourth case,
\[
\p_u^\beta H(0,0,0)=0\quad\mbox{for}\ \alpha+1 \leq \forall\beta \leq \beta_0<2\alpha
\quad\mbox{and}\  \d\int_{\R}g(x)dx=0,
 \]
 should be added in the result of the general theory (\ref{lifespan_general})
 with a new estimate,
 \begin{equation}
\label{lifespan_general_new}
\wt{T}(\e)\ge
 C\e^{-\beta_0(\alpha+1)/(\beta_0+2)}.
 \end{equation}
\end{theorem}

\par
We note that (\ref{lifespan_general_new}) is stronger than both the second and the third cases
in (\ref{lifespan_general}) because we have
\[
\frac{\beta_0(\alpha+1)}{\beta_0+2}>\frac{\alpha(\alpha+1)}{\alpha+2}
\quad\mbox{by}\ \alpha+1\le\beta_0
\]
and
\[
\frac{\beta_0(\alpha+1)}{\beta_0+2}>\frac{\beta_0}{2}
\quad\mbox{by}\ \beta_0<2\alpha.
\]
The sharpness of (\ref{lifespan_general_new}) can be verified by setting
\[
p=\alpha+1,\ q=0,\ r=\beta_0+1,
\]
or
\[
p+q=\alpha+1,\ r=\beta_0+1\ \mbox{with}\ q\neq0
\]
in Theorem \ref{thm:gcombined}.

\par
We shall omit the proof of Theorem \ref{thm:Takamatsu},
but point out that it is initiated by that of Thereom \ref{thm:gcombined},
especially the key is how to split $L^\infty$ norm of the solution itself
according to the appropriate domains.
Let us see it in the next section.

\subsection{Strategy of the proof of Theorem \ref{thm:gcombined}.} 
The key of our success to prove the generalized combined effect case,
especially the lower bound of the lifespan estimate in  (\ref{lifespan_gcombined}),
in Theorem \ref{thm:gcombined} is to handle
the system of integral equations of $(u-\e u^0,u_t-\e u^0_t)$, where
$u^0$ is a solution of the free wave equation with the same initial data $(u,u_t)|_{t=0}=(f,g)$;
\begin{equation}
\label{free}
u^0(x,t)=\frac{1}{2}\{f(x+t)+f(x-t)\}+\frac{1}{2}\int_{x-t}^{x+t}g(y)dy.
\end{equation}
The basic argument is due to John \cite{John79}
in which classical solutions of semilinear wave equations in three space dimensions
are constructed in a weighted $L^\infty$ space. 

\par
Assume that 
\begin{equation}
\label{supp_data}
\mbox{supp}\ f,g\subset\{|x|\le R\}\ \mbox{with}\ R>1,
\end{equation}
which yields that the finiteness of the propagation seed of the wave,
\begin{equation}
\label{supp_sol}
\mbox{supp\ }u(x,t)\subset\{|x|\le t+R\}
\end{equation}
by standard argument on small solutions of nonlinear wave equations.
For example, see Appendix in John \cite{John_book}.
Then, set an \lq\lq interior" domain $D$ by
\[
D:=\{t-|x|\ge R\}
\]
and define a sequence $\{(U_j,V_j)\}$ by
\[
\left\{
\begin{array}{l}
U_{j+1}=L(A_0|V_j+\e u_t^0|^p|U_j+\e u^0|^q+B_0|U_j+\e u^0|^r),\\
U_1=0,\\
V_{j+1}=L'(A_0|V_j+\e u_t^0|^p|U_j+\e u^0|^q+B_0|U_j+\e u^0|^r),\\
V_1=0,
\end{array}
\right.
\]
where $L(v)$ for a function $v$ of space-time variables is from Duhamel's term defined by 
\begin{equation}
\label{L}
\begin{array}{l}
\d L(v)(x,t):=\frac{1}{2}\int_0^tds\int_{x-t+s}^{x+t-s}v(y,s)dy,\\
\d L'(v)(x,t):=\frac{\p}{\p t}L(v)(x,t).
\end{array}
\end{equation}
This sequence $\{(U_j,V_j)\}$ will converge to $(u-\e u^0,u_t-\e u_t^0)$
in a function space
\[
\begin{array}{ll}
X:=&\{(U,V)\in\{C^1(\R\times[0,T])\}^2\\
&\quad:\|(U,V)\|_X<\infty,\ \mbox{supp}\ (U,V)\subset\{|x|\le t+R\}\},
\end{array}
\]
where
\[
\begin{array}{l}
\|(U,V)\|_X:=\|U\|_1+\|U_x\|_1+\|V\|_2+\|V_x\|_2,\\
\|U\|_1:=\d\sup_{(x,t)\in\R\times[0,T]}(t+|x|+R)^{-1}|U(x,t)|,\\
\|V\|_2:=\d\sup_{(x,t)\in\R\times[0,T]}\{\chi_D(x,t)\\
\qquad+(1-\chi_D(x,t))(t+|x|+R)^{-1}\}|V(x,t)|
\end{array}
\]
and $\chi_D$ is a characteristic function of the interior domain $D$.
In view of the expression of $u^0$ in (\ref{free}),
it follows from (\ref{supp_data}) and the assumption on $g$ yield the strong Huygens' principle,
\begin{equation}
\label{Huygens}
u^0(x,t)\equiv\int_{\R}g(x)dx=0\quad\mbox{in}\ D.
\end{equation}
\par
Then we have main a priori estimates;
\[
\begin{array}{ll}
\|L(|V|^p|U|^q)\|_1\le C\|V\|_2^p\|U\|_1^q(T+R)^{p+q},\\
\|L(|U|^r)\|_1\le C\|U\|_1^r(T+R)^{r+1},\\
\|L'(|V|^p|U|^q)\|_2\le C\|V\|_2^p\|U\|_1^q(T+R)^{p+q},\\
\|L'(U|^r)\|_2\le C\|U\|_1^r(T+R)^{r+1}
\end{array}
\]
with some positive constant $C$ independent of $\e$ and $T$
which give us the key estimates for $p+q<r$;
\[
\begin{array}{ll}
\|U_{j+1}\|_1\lesssim &\e^{p+q}+\|V_j\|_2^p\|U_j\|_1^qT^{p+q}
+\|U_j\|_1^rT^{r+1}+\mbox{harmless terms},\\
\|V_{j+1}\|_2 \lesssim & \mbox{the same as above},
\end{array}
\]
where $\lesssim$ stands for mod constant independent of $\e$ and $T$.
Therefore the boundedness of the sequence, $\|U_j\|_1,\|V_j\|_2\lesssim\e^{p+q}$ for all $j\in\N$,
follows from the competition between the second and third terms in the right-hand side.
When $(r+1)/2<p+q$, the third term is major, so that
\[
\e^{r(p+q)}T^{r+1}\lesssim \e^{p+q}
\]
is the required condition.
The difference to convergence of the sequence as well as those for the spatial derivative
is almost the same as the boundedness above.
Therefore we obtain the desired lower bound;
\[
T(\e)\ge C\e^{-(p+q)(r-1)/(r+1)}.
\]

\par
 For the upper bound in this generalized combined effect,
 its proof is easy if we follow the argument by higher dimensional case by
 Han and Zhou \cite{HZ14}.
 In fact, assume that $f(x)\ge0(\not\equiv0),\ g(x)\equiv0$.
Set
\[
F(t):=\d\int_{\R}u(x,t)dx.
\]
Then, the equation and (\ref{supp_sol}) imply that
\[
F''(t)=\int_{\R}(A_0|u_t|^p|u|^q+B_0|u|^r)dx.
\]
The second term and H\"older's inequality yield that
\begin{equation}
\label{frame}
F''(t)\gtrsim (t+R)^{-(r-1)}|F(x)|^r
\quad\mbox{for}\ t\ge0.
\end{equation}
This is the main inequality.
Next we shall make use of the first term to obtain the better estimate than the single term.
For this purpose, assume further that
\[
f(x),-f'(x)\ge \exists c>0\ \quad\mbox{for}\ x\in(-R/2,0).
\]
Then, we have that
\[
\left\{
\begin{array}{l}
u(x,t)\gtrsim\e f(x-t)\gtrsim\e,\\
u_t(x,t)\gtrsim-\e f'(x-t)\gtrsim\e
\end{array}
\right.
\mbox{for}\ t<x<t+R.
\]
Plugging these estimates into the first term in (\ref{frame}), we obtain that
\[
F''(t)\gtrsim \e^{p+q}\quad\mbox{for large $t$}
\]
which implies that
\begin{equation}
\label{first_step}
F(t)\gtrsim \e^{p+q}t^2\quad\mbox{for large $t$}.
\end{equation}
These estimates (\ref{frame}) and (\ref{first_step}) will give us the desired result
according to the improved Kato's lemma by Takamura \cite{Takamura15}.
But here we show a brief feeling that the result is correct as follows.
Plugging (\ref{first_step}) into (\ref{frame}), we obtain that
\[
\begin{array}{l}
F''(t)\gtrsim\e^{r(p+q)}t^{r+1}\quad\mbox{for large $t$}
\end{array}
\]
which implies that
\begin{equation}
\label{second_step}
\begin{array}{l}
F(t)\gtrsim \e^{r(p+q)}t^{r+3}\quad\mbox{for large $t$}.
\end{array}
\end{equation}
The difference of the lower bound of $F$ between (\ref{first_step}) and (\ref{second_step}) is 
\[
\e^{(p+q)(r-1)}t^{r+1}.
\]
This quantity is larger than some constant, we will reach to the desired lifespan estimate
\[
T(\e)\le C\e^{-(p+q)(r-1)/(r+1)}
\]
by standard iteration argument.
Because we don't have to cut the time interval in each step,
which means that this procedure can be repeated infinitely many times
in the fixed time interval.
In this way, we obtain the sharp estimate of the lifespan in the generalized combined effect case. 

\par
Other cases except for the generalized combined effect are similar to those above.
However, we point out that the existence part under the assumption of
$\d\int_{\R}g(x)dx\neq0$ is most difficult practically among all the proofs
in choosing appropriate weight functions.

\section{The case of variable coefficients} 
In this section, we assume for our model (\ref{IVP_gcombined}) that
$A(x,t)$, or $B(x,t)$, is of
\[
\frac{1}{\langle t+\langle x\rangle\rangle^{1+a}
\langle t-\langle x\rangle\rangle^{1+b}\langle x\rangle^{1+c}},\ \mbox{or}\ \equiv0,
\]
where $a,b,c\in\R$ and $\langle x\rangle:=\sqrt{1+x^2}$.
The reason to take many $\langle\ \rangle$s is to ensure the differentiability of
$A,B$ to construct a classical solution.
We will see that this kind of models will be a key in extending the general theory
for (\ref{GIVP}) to the non-autonomous terms as stated at the end of Introduction.
The weights of this kind were firstly introduced by 
Belchev, Kepka and Zhou \cite{BKZ2002}, later Liu and Zhou \cite{LZ2008}
in higher space dimensions to show a blow-up result.
But they set a special weight to make use of some geometric transform
to absorb it, and reduced the equation to ordinary differential inequality of the functional
without any argument on the local existence of the solution as well as the lifespan estimates.

\subsection{Motivation of the problem}
Our motivation to consider such $A$ and $B$ above
comes from an initial value problem for the semilinear damped wave equations;
\begin{equation}
\label{IVP_damped}
\left\{
\begin{array}{ll}
\d v_{tt}- v_{xx}+\frac{2}{1+t} v_t=|v|^p &\quad \mbox{in}\ \R\times (0,\infty),\\
 v(x,0)= \e f(x) , \quad v_t(x,0)= \e g(x) &\quad \mbox{for}\ x\in\R
\end{array}
\right.
\end{equation}
with the same setting on the initial data.
This problem is a very important model as it has a critical decay in damping term,
namely scaling invariant damping.
The time-decay $(1+t)^{-1}$ is a threshold between heat-like
with weaker decay and wave-like with stronger decay
in the sense that the each critical exponent is the same as
Fujita one and Strauss one respectively.
In the critical case, the size of the constant in front of the damping term is important.
There is also a threshold on the constant and this case \lq\lq2" is in the domain of heat-like
in one dimension.
See the introduction of Kato, Takamura and Wakasa \cite{KTW19}
for precise results and references.
In fact, D'Abbicco \cite{D'Abbico15} showed for $p>3$
that the energy solution of (\ref{IVP_damped}) exists globally-in-time,
while Wakasugi \cite{Wakasugi14} obtained its counter part of
finite-time blow-up of the energy solution for $1<p\le3$. 
This critical exponent $3$ is Fujita one in one space dimension,
so one may expect that the solution behaves like one of semilinear heat equations
for which $u_{tt}$ is neglected from (\ref{IVP_damped}).
But, this is not true.

\par
In fact, Liouville transform $u(x,t)=(1+t)v(x,t)$ shows that (\ref{IVP_damped}) is equivalent to
\begin{equation}
\label{IVP_weighted}
\left\{
\begin{array}{ll}
\d u_{tt}- u_{xx}=\frac{|u|^p}{(1+t)^{p-1}} & \mbox{in}\ \R\times (0,\infty),\\
 u(x,0)= \e f(x),\ u_t(x,0)= \e\{f(x)+g(x)\} & \mbox{for}\ x\in\R,
\end{array}
\right.
\end{equation}
so that all the technique for semilinear wave equations are applicable to this problem,
and we have the following results on the lifepsan estimates for (\ref{IVP_damped}).
Wakasa \cite{Wakasa16} obtained that
\begin{equation}
\label{wak-est}
\begin{array}{c}
T(\e)\sim
\left\{
\begin{array}{lll}
C\e^{-(p-1)/(3-p)} &\mbox{for}&\d 1<p<3,\\
\exp(C\e^{-(p-1)}) &\mbox{for}&\d p=3 \\
\end{array}
\right. 
\\
\mbox{if}\ \d \int_{\R}\{f(x)+g(x)\}dx \neq0
\end{array}
\end{equation}
This is the heat-like estimate in the sense that 
it coincides with the corresponding semilinear heat equations 
for which $u_{xx}$ is neglected in (\ref{IVP_damped}).  
In the original paper \cite{Wakasa16}, the condition on the initial data is missing,
but it was expected to be natural as the critical exponent is Fujita one.
Later, Kato, Takamura and Wakasa \cite{KTW19} proved that
\begin{equation}
\label{KTW-est}
\begin{array}{c}
T(\e)\sim
\left\{
\begin{array}{lll}
C\e^{-p(p-1)/(1+2p-p^2)} &\mbox{for}&\d 1<p<2,\\
Cb(\e) &\mbox{for}&\d p=2,\\
C\e^{-p(p-1)/(3-p)} &\mbox{for}&\d 2<p<3,\\
\exp(C\e^{-p(p-1)}) &\mbox{for}&\d p=3
\end{array}\right. \\
\mbox{if}\ \d \int_{\R}\{f(x)+g(x)\}dx =0,
\end{array}
\end{equation}
where $b=b(\e)$ is a positive number satisfying $\e^2 b \log(1+b)=1$. 
This is the wave-like estimate.
In deed, the cases $1<p<2$ and $p=3$ are
the same forms as those of 3-dimensional semillinear wave equations. 

\par
In this way, it is very important to study the weighted nonlinear terms,
and it is natural to extend the results to more general weighted terms than (\ref{IVP_weighted}).
The first simple question may go to the case of $x$-decay.

\subsection{Spatially weighted nonlinear terms.}
First we consider the following case of spatial weights;
\begin{equation}
\label{x-depend_u}
A(x,t)\equiv0\quad\mbox{and}\quad B(x,t)=\frac{1}{\langle x\rangle^{1+c}},
\end{equation}
where $c\in\R$, in (\ref{IVP_gcombined}).
This setting was first introduced by Suzuki \cite{Suzuki10} under the supervision
by Prof. M. Ohta (Tokyo Univ. of Sci. Japan),
but with the assumption that the initial data is of non-compact support.
Later, Kubo, Osaka and Yazici \cite{KOY13} improved the result and
Wakasa \cite{Wakasa17} finalized it.
For the compactly supported case, we have the following result.

\begin{theorem}[Kitamura, Morisawa and Takamura \cite{KMT22}]
\label{thm:KMT22}
Assume (\ref{x-depend_u}).
Then, the lifespan of a classical solution of (\ref{IVP_gcombined}) satisfies the following estimates.
\begin{equation}
\label{lifespan_non-zeroKMT}
T(\e)\sim
\left\{
\begin{array}{ll}
C\e^{-(p-1)/(1-c)} & \mbox{for}\ c<0,\\
\phi^{-1}(C\e^{-(p-1)}) & \mbox{for}\ c=0,\\
C\e^{-(p-1)} & \mbox{for}\ c>0
\end{array}
\right.
\quad
\mbox{if}\ \int_\R g(x)dx\neq0
\end{equation}
and
\begin{equation}
\label{lifespan_zeroKMT}
T(\e)\sim
\left\{
\begin{array}{ll}
C\e^{-p(p-1)/(1-pc)} & \mbox{for}\ c<0,\\
\psi^{-1}(C\e^{-p(p-1)}) & \mbox{for}\ c=0,\\
C\e^{-p(p-1)} & \mbox{for}\ c>0
\end{array}
\right.
\quad
\mbox{if}\ \int_\R g(x)dx=0,
\end{equation}
where $\phi^{-1}$ and $\psi^{-1}$ are inverse functions defined by $\phi(s)=s\log(2+s)$ and $\psi(s)=s\log^p(2+s)$, respectively. 
\end{theorem}

\begin{remark}
Wakasa \cite{Wakasa17} established the estimate in (\ref{lifespan_non-zeroKMT})
with the assumption that $f\in L^\infty(\R),\ g\in L^1(\R)$ and $c\ge-1$.
The last condition ensures the existence of local-in-time solutions
with non-compactly supported data.
We also remark that $|u|^p$ in the nonlinear term can be replaced with
$|u|^{p-1}u$ in \cite{Wakasa17} due to the fact that
the positiveness of the solution can be obtained easier than the compactly supported case.
\end{remark}

\par
The strategy of the proof of Theorem \ref{thm:KMT22} for the existence part is similar to the one
of Theorem \ref{thm:gcombined} based on weighted $L^\infty$ estimates of the solution,
so we don't comment about it here.
For the blow-up part, the iteration argument of the point-wise estimate is employed.
The functional method like the proof of Theorem \ref{thm:gcombined} cannot be applied
to this case due to the effect of the weight.

\par
After the work \cite{KMT22}, one may study the counter case of
\begin{equation}
\label{x-depend_u_t}
A(x,t)=\frac{1}{\langle x\rangle^{1+c}}\quad\mbox{and}\quad B(x,t)\equiv0,
\end{equation}
where $c\in\R$, in (\ref{IVP_gcombined}) with $q=0$.
For this equation, one may expect that the result is not so interesting as
$|x|\sim t$ in dealing with the nonlinear term $|u_t|^p$.
But the result is 

\begin{theorem}[Zhou \cite{Zhou01}, Kitamura, Morisawa and Takamura \cite{KMT23}]
\label{thm:KMT23}
Assume (\ref{x-depend_u_t}) with $q=0$.
Then, the lifespan of a classical solution of (\ref{IVP_gcombined}) satisfies the following estimates.
\begin{equation}
\label{lifespan_KMT23}
\begin{array}{l}
T(\e)\sim
\left\{
\begin{array}{ll}
C\e^{-(p-1)/(-c)} & \mbox{for}\ c<0,\\
\exp(C\e^{-(p-1)}) & \mbox{for}\ c=0,
\end{array}
\right.
\\
T(\e)=\infty\quad\mbox{for}\ c>0.
\end{array}
\end{equation}
\end{theorem}

\par
We note that Zhou \cite{Zhou01} established the blow-up part with $c=-1$,
namely non-weighted case.
The proof of Theorem \ref{thm:KMT23} is similar to the one
of Theorem \ref{thm:KMT22}, so we don't comment about it here also.
The main difference between (\ref{x-depend_u}) and (\ref{x-depend_u_t})
is a possibility to obtain the global-in-time existence for (\ref{x-depend_u_t})
while there is no such a situation due to the effect of the origin in space for (\ref{x-depend_u}).

\par
According to the result on (\ref{IVP_damped})
and Theorems \ref{thm:KMT22} and \ref{thm:KMT23},
one may expect that it is sufficient to study the weights in the nonlinear terms
by powers of $(1+t)$ or $\langle x\rangle$ for the purpose to extend the general theory,
but we have to take into account of \lq\lq characteristic weights"
in handling nonlinear wave equations.
Let's see it in the next subsection.

\subsection{Weighted nonlinear terms in the characteristic directions.}
\par
It is well-known that the wave propagates along with the characteristic directions,
so that it is not sufficient to study the weights of powers by $(1+t)$ or $\langle x\rangle$ only.
Therefore, as a breakthrough to this direction, we set

\begin{equation}
\label{xt-depend_u}
A(x,t)\equiv0\quad\mbox{and}\quad B(x,t)=\frac{1}{\langle t+\langle x\rangle\rangle^{1+a}
\langle t-\langle x\rangle\rangle^{1+b}},
\end{equation}
where $a,b\in\R$, in (\ref{IVP_gcombined}).
Then, we obtained the following result.

\begin{theorem}[Kitamura, Wakasa and Takamura \cite{KTW23}]
\label{thm:KTW23}
Assume (\ref{xt-depend_u}).
Then, the lifespan $T(\e)$ of the classical solution of (\ref{IVP_gcombined}) satisfies the following estimates;
\begin{equation}
\label{main_global}
T(\e)=\infty\quad \mbox{for}\ a+b>0\ \mbox{and} \ a>0,
\end{equation}
and
\begin{equation}
\label{main_lifespan1}
T(\e)\sim
\left\{
\begin{array}{lllll}
\exp(C\e^{-(p-1)}) &\mbox{for}&
\begin{array}{l}
a+b=0\ \mbox{and}\ a>0,\\
\mbox{or}\ a=0\ \mbox{and} \ b>0,
\end{array}
\\ 
\exp(C\e^{-(p-1)/2}) &\mbox{for}& a=b=0, \\
C\e^{-(p-1)/(-a)} &\mbox{for}& a<0 \ \mbox{and}\ b>0,\\
\phi_1^{-1}(C\e^{-(p-1)})&\mbox{for}& a<0\ \mbox{and}\ b=0,\\
C\e^{-(p-1)/(-a-b)}&\mbox{for}& a+b<0\ \mbox{and}\ b<0\\
\end{array}
\right.
\end{equation}
if 
\[
\int_{\R}g(x)dx\neq0,
\]
where $\phi_1^{-1}$ is an inverse function defined by
\begin{equation}
\label{phi1}
\phi_1(s)=s^{-a}\log(2+s).
\end{equation}
On the other hand, it holds that
\begin{equation}
\label{main_lifespan2}
T(\e)\sim
\left\{
\begin{array}{lllll}
\exp(C\e^{-(p-1)}) &\mbox{for}&  a=0\ \mbox{and}\ b>0,\\
\exp(C\e^{-p(p-1)}) &\mbox{for}& a+b=0\ \mbox{and} \ a>0, \\  
\exp(C\e^{-p(p-1)/(p+1)}) &\mbox{for}& a=b=0, \\
C\e^{-(p-1)/(-a)} &\mbox{for}& a<0 \ \mbox{and}\ b>0,\\
\psi_1^{-1}(C\e^{-p(p-1)})&\mbox{for}& a<0\ \mbox{and}\ b=0,\\
C\e^{-p(p-1)/(-pa-b)} &\mbox{for}& a<0\ \mbox{and}\ b<0,\\
\psi_2^{-1}(C\e^{-p(p-1)})&\mbox{for}& a=0\ \mbox{and}\ b<0,\\
C\e^{-p(p-1)/(-a-b)} &\mbox{for}& a+b<0 \ \mbox{and}\ a>0
\end{array}
\right.
\end{equation}
if
\[
\int_{\R}g(x)dx=0,
\]
where $\psi_1^{-1}$ and $\psi_2^{-1}$ are inverse functions defined by 
\begin{equation}
\label{psi1_psi2}
\psi_1(s)=s^{-pa}\log(2+s)\ \mbox{and} \ \psi_2(s)=s^{-b}\log^{p-1}(2+s).
\end{equation}
\end{theorem}

\begin{remark}
The estimates (\ref{main_lifespan2}) with $a=p-2$ and $b=-1$ coincide
with (\ref{wak-est}) and (\ref{KTW-est})
because $(1+t)$ is equivalent to $\LR{t+\LR{x}}$ 
by finite propagation speed of the wave like (\ref{supp_sol}).
\end{remark}

\par
The strategy of the proof of Theorem \ref{thm:KTW23} is also almost the same as the one
of Theorem \ref{thm:KMT22}, but the weights cause many technical difficulties,
but we shall skip the details fot the purpose of this paper.
The main concern of Theorem \ref{thm:KTW23} is that
we have interactions among two characteristic directions as the critical line
$a+b=0$ with $b\le0$ and $a=0$ with $b\ge0$ which divides $ab$-plane into
two domains of the global-in-time existence and the blow-up in finite time.

\par
In contrast, if

\begin{equation}
\label{xt-depend_u_t}
A(x,t)=\frac{1}{\langle t+\langle x\rangle\rangle^{1+a}
\langle t-\langle x\rangle\rangle^{1+b}}\quad\mbox{and}\quad B(x,t)\equiv0,
\end{equation}
where $a,b\in\R$, in (\ref{IVP_gcombined}) with $q=0$,
then we obtained the different critical line on $ab$-plane
from Theorem \ref{thm:KTW23} as follows.

\begin{theorem}[Kitamura \cite{Kitamura}]
\label{thm:Kitamura}
Assume (\ref{x-depend_u_t}) with $q=0$.
Then, the lifespan of a classical solution of (\ref{IVP_gcombined}) satisfies the following estimates.
\begin{equation}
\label{lifespan_Kitamura}
\begin{array}{l}
T(\e)\sim
\left\{
\begin{array}{l}
C\e^{-(p-1)/(-a)}\\
\qquad \mbox{for}\ a<0\ \mbox{and}\ b\ge-p,\\
C\e^{-p(p-1)/\{-p(1+a)-b\}}\\
\qquad\mbox{for}\ p(1+a)+b<0\ \mbox{and}\ b<-p,\\
\exp(C\e^{-(p-1)})\\
\qquad\mbox{for}\ a=0\ \mbox{and}\ b\ge-p,\\
\exp(C\e^{p(p-1)})\\
\qquad\mbox{for}\ a>0\ \mbox{and}\ p(1+a)+b=0,
\end{array}
\right.
\\
T(\e)=\infty\quad\mbox{for}\ a>0\ \mbox{and}\ p(1+a)+b>0.
\end{array}
\end{equation}
\end{theorem}

\par
We note again that all the estimates above are established
whatever the value of $\d\int_{\R}g(x)dx$ is, due to the nonlinear term $|u_t|^p$.
The strategy of the proof of Theorem \ref{thm:Kitamura} is also almost the same as the one
of Theorem \ref{thm:KMT22}.
\vskip10pt
\par
{\bf Concluding remark}\quad
The series of our theorems above are major
if one intends to extend the general theory for non-autonomous equations.
But they are not sufficient still as we never try to analyze the combined effect
for variable coefficients case including the different order of $A$ and $B$.

\subsection{Other models.}
\par
Finally we shall comment on the different nonlinear terms.
Recently we obtained the following result.

\begin{theorem} [Sasaki, Takamatsu and Takamura \cite{STT23}]
\label{thm:STT23}
Assume that $A\equiv1,B\equiv0,q=0$ in (\ref{IVP_gcombined}).
Moreover, $|u_t|^p$ is replaced with $|u_x|^p$.
Then we have that
\begin{equation}
\label{lifespan_STT23}
T(\e)\sim C\e^{-(p-1)}.
\end{equation}
\end{theorem}

\par
One may feels that this result is trivial according to Theorem \ref{thm:KMT22} with $c=-1$.
But the proof is different from each others.
Especially there is a difficulty on the blow-up part 
as the point-wise positivity of $u_x$ cannot be obtained by
\[
\frac{\p}{\p x}L(v)(x,t),
\]
where $L(v)$ is the one in (\ref{L}). 
This situation is overcome by taking suitably weighted functional of the solution
which was introduced by Rammaha \cite{Rammaha95, Rammaha97}.
We believe that Theorem \ref{thm:STT23} will contribute to analysis
of the \lq\lq blow-up boundary" for the equation
\[
u_{tt}-u_{xx}=|u_x|^p
\]
for which there is no result til now.
See \cite{STT23} for references therein to the blow-up boundary.



\begin{thebibliography}{99}


\bibitem{D'Abbico15}{M. D'Abbicco},
{\it The threshold of effective damping for semilinear wave equations},
Mathematical Methods in Applied Sciences, {\bf 38} (2015), 1032-1045.

\bibitem{BKZ2002}{E. Belchev, M. Kepka and Z. Zhou},
{\it Finite-time blow-up of solutions to semilinear wave equations},
Special issue dedicated to the memory of I. E. Segal. J. Funct. Anal. {\bf 190} (2002), no. 1, 233-254.

\bibitem{BS}{R. Bianchini and G. Staffilani},
{\it Revisitation of Tartar's result on a semilinear hyperbolic system with null condition},
arXiv: 2001.03688v2 [math.AP].

\bibitem{DFW19}{W. Dai, D. Fang and C. Wang},
{\it Global existence and lifespan for semilinear wave equations with mixed nonlinear terms},
J. Differential Equations, {\bf 267} (2019), no. 5, 3228-3354.

\bibitem{HZ14}{W. Han and Y. Zhou},
{\it Blow up for some semilinear wave equations in multi-space dimensions},
Comm. Partial Differential Equations, {\bf 39} (2014), no. 4, 651-665.

\bibitem{HWY16}{K. Hidano,  C. Wang and K. Yokoyama},
{\it Combined effects of two nonlinearities in lifespan
of small solutions to semi-linear wave equations},
Math.  Ann., {\bf 366} (2016), no. 1-2, 667-694.

\bibitem{IKTW}{T. Imai, M. Kato, H. Takamura and K. Wakasa}, 
{\it The sharp lower bound of the lifespan of solutions to semilinear wave equations
with low powers in two space dimensions},
Asymptotic analysis for nonlinear dispersive and wave equations, 31-53,
Adv. Stud. Pure Math., {\bf 81}, Math. Soc. Japan, Tokyo, 2019. 

\bibitem{John79}{F. John},
{\it Blow-up of solutions of nonlinear wave equations in three space dimensions},
Manuscripta Math., {\bf 28} (1979), no. 1-3, 235-268.

\bibitem{John_book}{F. John},
\lq\lq Nonlinear Wave Equations, Formation of Singularities",
ULS Pitcher Lectures in Mathematical Science, Lehigh University,
American Mathematical Society, Providence, RI, 1990.

\bibitem{Katayama01}{S. Katayama},
{\it Lifespan of solutions for two space dimensional wave equations
with cubic nonlinearity},
Comm. Partial Differential Equations, {\bf 26} (2001), no. 1-2, 205-232.

\bibitem{KTW19}{M. Kato, H. Takamura and K. Wakasa},
{\it The lifespan of solutions of semilinear wave equations
with the scale-invariant damping in one space dimension},
Differential Integral Equations, {\bf 32} (2019), no. 11-12, 659-678.

\bibitem{KSTT}{R. Kido, T. Sasaki, S. Takamatsu and H. Takamura},
{\it The generalized combined effect for one dimensional wave equations
with semilinear terms including product type},
arXiv:2305.00180v2 [math.AP].

\bibitem{Kitamura}{S. Kitamura},
{\it Semilinear wave equations of derivative type with characteristic weights
in one space dimension},
arXiv:2211.12295v2 [math.AP].

\bibitem{KMT22}{S. Kitamura, K. Morisawa and H. Takamura},
{\it The lifespan of classical solutions of semilinear wave equations
with spatial weights and compactly supported data in one space dimension},
J. Differential Equations, {\bf 307} (2022), 486-516.

\bibitem{KMT23}{S. Kitamura, K. Morisawa and H. Takamura},
{\it Semilinear wave equations of derivative type
with spatial weights in one space dimension},
Nonlinear Analysis, RWA, {\bf 72} (2023), Paper No.103764.

\bibitem{KTW23}{S. Kitamura, H. Takamura and K. Wakasa},
{\it The lifespan estimates of classical solutions of one dimensional semilinear wave equations
with characteristic weights},
J. Math. Anal. Appl. {\bf528} (2023), Paper No. 127516.

\bibitem{KOY13}{H. Kubo, A.Osaka and M.Yazici},
{\it Global existence and blow-up for wave equations
with weighted nonlinear terms in one space dimension},
Interdisciplinary Information Sciences, {\bf 19} (2013), 143-148.

\bibitem{LZ2008}{X. Liu and Y. Zhou},
{\it Global nonexistence of solutions to a semilinear wave equation in the Minkowski space},
Appl. Math. Lett. {\bf 21} (2008), no. 8, 849-854.

\bibitem{LYZ91}{T.-T. Li (D.-Q. Li), X. Yu and Y. Zhou},
{\it Dur\'ee de vie des solutions r\'eguli\`eres pour les \'equations des ondes non lin\'eaires unidimensionnelles} (French),
C. R. Acad. Sci. Paris S\'er. I Math., {\bf 312} (1991), no. 1, 103-105.

\bibitem{LYZ92}{T.-T. Li, X. Yu and Y. Zhou},
{\it Life-span of classical solutions to one-dimensional nonlinear wave equations},
Chinese Ann. Math., Ser. B, {\bf13} (1992), no. 3, 266-279. 

\bibitem{LW20}{M. Liu and C. Wang},
{\it Blow up for small-amplitude semilinear wave equations
with mixed nonlinearities on asymptotically Euclidean manifolds},
J. Differential Equations, {\bf269} (2020), no. 10, 8573-8596. 

\bibitem{LYY18}{G. K. Luli, S. Yang, P. Yu},
{\it On one-dimension semi-linear wave equations with null conditions},
Adv. Math., {\bf 329} (2018), 174-188.

\bibitem{Nakamura14}{M. Nakamura},
{\it Remarks on a weighted energy estimate and its application
to nonlinear wave equations in one space dimension},
J. Differential Equations, {\bf 256}, no. 2, (2014), 389-406.

\bibitem{MST}{K. Morisawa, T. Sasaki and H. Takamura},
{\it The combined effect in one space dimension beyond the general theory
for nonlinear wave equations},
Comm. Pure Appl. Anal., {\bf 22} (2023), no. 5, 1629-1658.

\bibitem{MST_erratum}{K. Morisawa, T. Sasaki and H. Takamura},
{\it Erratum to \lq\lq The combined effect in one space dimension beyond the general theory
for nonlinear wave equations"}, Comm. Pure Appl. Anal., {\bf22} (2023), no. 10, 3200-3202.

\bibitem{Rammaha95}{M. A. Rammaha},
{\it Upper bounds for the life span of solutions to systems of nonlinear wave equations
in two and three space dimensions},
Nonlinear Anal. {\bf 25} (1995), no. 6, 639-654.

\bibitem{Rammaha97}{M. A. Rammaha},
{\it A note on a nonlinear wave equation in two and three space dimensions},
Comm. Partial Differential Equations {\bf 22} (1997), no. 5-6, 799-810. 

\bibitem{STT23}{T. Sasaki S. Takamatsu and H. Takamura},
{\it The lifespan of classical solutions of one dimensional wave equations
with semilinear terms of the spatial derivative},
AIMS Math., {\bf 8}, no. 11, 25477-25486.

\bibitem{Suzuki10}{A. Suzuki},
\lq\lq Global Existence and Blow-Up of solutions
to Nonlinear Wave Equations in One Space Dimension" (in Japanese),
Master Thesis, Saitama University, 2010.

\bibitem{Takamatsu}{S. Takamatsu},
{\it Improvement of the general theory for one dimensional nonlinear wave equations
related to the combined effect}, arXiv: 2308.02174v2 [math.AP].  

\bibitem{Takamura15}{H. Takamura},
{\it Improved Kato's lemma on ordinary differential inequality and its applications
to semilinear wave equations},
Nonlinear Anal.,  {\bf 125} (2015), 227-240.

\bibitem{Tartar}{L. Tartar},
{\it Some existence theorems
for semilinear hyperbolic systems in one space variable},
1981, unpublished.

\bibitem{Wakasa16}{K. Wakasa},
{\it The lifespan of solutions to semilinear damped wave equations in one space dimension},
Communications on Pure and Applied Analysis, {\bf 15} (2016), 1265-1283.

\bibitem{Wakasa17}{K. Wakasa},
{\it The lifespan of solutions to wave equations with weighted nonlinear terms in one space dimension},
Hokkaido Math. J., {\bf 46} (2017), 257-276.

\bibitem{Wakasugi14}{Y. Wakasugi},
{\it Critical exponent for the semilinear wave equation with scale invariant damping}, Fourier analysis, 375-390, Trends Math., Birkh\"{a}user/Springer, Cham, 2014.

\bibitem{Zha20}{D. Zha},
{\it On one-dimension quasilinear wave equations with null conditions},
Calc. Var. Partial Differential Equations,  {\bf 59} (2020), no. 3, Paper No. 94, 19 pp.

\bibitem{Zha22}{D. Zha},
{\it Global stability of solutions to two-dimension and one-dimension systems
of semilinear wave equations},
J. Funct. Anal., {\bf 282} (2022), no. 1, Paper No. 1092219, 26 pp.

\bibitem{Zhou92}{Y. Zhou},
{\it Life span of classical solutions to $u_{tt}-u_{xx}=|u|^{1+\alpha}$},
Chinese Ann. Math. Ser.B, {\bf 13} (1992), no. 2, 230-243.

\bibitem{Zhou01}{Y. Zhou},
{\it Blow up of solutions to the Cauchy problem for nonlinear wave equations},
Chinese Ann. Math. Ser. B, {\bf22} (2001), no. 3, 275-280.

\bibitem{Zhou}{Y. Zhou},
{\it Blow up of solutions to the Cauchy problem for nonlinear wave equations},
preprint, Fudan University, 1992.


\end{thebibliography}
\end{document}